\documentclass{article}%
\usepackage{amssymb}%
\usepackage{amsmath}%
\usepackage{amsfonts}%
\usepackage{graphicx}
\newtheorem{theorem}{Theorem} [section]

\newtheorem{corollary}[theorem]{Corollary}

\newtheorem{lemma}[theorem]{Lemma}

\newtheorem{proposition}[theorem]{Proposition}

\newenvironment{proof}[1][Proof]{\noindent\textbf{#1.} }{\ \rule{0.5em}{0.5em}}
\begin{document}

\title{The Intersection of All Maximum Stable Sets of a Tree and its Pendant Vertices}
\author{Vadim E. Levit and Eugen Mandrescu\\Department of Computer Science\\Holon Academic Institute of Technology\\52 Golomb Str., P.O. Box 305\\Holon 58102, ISRAEL\\\{levitv, eugen\_m\}@hait.ac.il}
\maketitle

\begin{abstract}
A \textit{stable} set in a graph $G$ is a set of mutually non-adjacent
vertices, $\alpha(G)$ is the size of a maximum stable set of $G$, and
$core(G)$ is the intersection of all its maximum stable sets. In this paper we
demonstrate that in a tree $T$, of order $n\geq2$, any stable set of size
$\geq n/2$ contains at least one pendant vertex. Hence, we deduce that any
maximum stable set in a tree contains at least one pendant vertex. We give a
new proof for a theorem of Hopkins and Staton \cite{HopStat} characterizing
strong unique trees. Using this result we show that if $\{A,B\}$ is the
bipartition of a tree $T$ and $S$ is a stable set with $\left\vert
S\right\vert >\min\{\left\vert A\right\vert ,\left\vert B\right\vert \}$, then
$S$ contains at least a pendant vertex.

Our main finding is the theorem claiming that if $T$ is a tree of order
$n\geq2$ that does not own a perfect matching (i.e., $2\alpha(T)>n$), then at
least two pendant vertices an even distance apart belong to $core(T)$. While
it is known that if $G$ is a connected bipartite graph of order $n$ $\geq2$,
then $\left\vert core(G)\right\vert \neq1$ (see Levit, Mandrescu
\cite{LeviMan2}), our new statement reveals an additional structure of the
intersection of all maximum stable sets of a tree. The above assertions give
refining of one assertion of Hammer, Hansen and Simeone \cite{HamHanSim}
stating that if a graph $G$ is of order less than $2\alpha(G)$, then $core(G)$
is non-empty, and also of a result of Jamison \cite{Jamison}, Gunter, Hartnel
and Rall \cite{GunHarRall}, and Zito \cite{Zito}, saying that for a tree $T$
of order at least two, $\left\vert core(T)\right\vert \neq1$.

\end{abstract}

\section{Introduction}

Throughout this paper $G=(V,E)$ is a simple (i.e., a finite, undirected,
loopless and without multiple edges) graph with vertex set $V=V(G)$, edge set
$E=E(G)$, and its order is $\left\vert V\right\vert $. If $X\subset V$, then
$G[X]$ is the subgraph of $G$ spanned by $X$. By $G-W$ we mean the subgraph
$G[V-W]$, if $W\subset V(G)$. We also denote by $G-F$ the partial subgraph of
$G$ obtained by deleting the edges of $F$, for $F\subset E(G)$, and we use
$G-e$, if $F$ $=\{e\}$. If $X,Y$ are non-empty disjoint subsets of $V$, then
$(X,Y)$ equals $\{xy\in E:x\in X,y\in Y\}$. Let $K_{n},P_{n}$ denote the
complete graph on $n\geq1$ vertices, and the chordless path on $n\geq2$ vertices.

A set $A\subseteq V$ is \textit{stable} if no two vertices from $A$ are
adjacent. A stable set of maximum size will be referred as to a
\textit{\ maximum stable set} of $G$, and the \textit{stability number }of
$G$, denoted by $\alpha(G)$, is the cardinality of a maximum stable set in
$G$. Let $\Omega(G)$ stand for the set $\{S:S$ \textit{is a maximum stable set
of} $G\}$, $core(G)=\cap\{S:S\in\Omega(G)\}$, and $\xi\left(  G\right)
=\left\vert core(G)\right\vert $ (see \cite{LeviMan3}).

The \textit{neighborhood} of a vertex $v\in V$ is the set $N(v)=\{w:w\in V$
\ \textit{and} $vw\in E\}$, while the \textit{closed neighborhood} of $v\in V
$ is $N[v]=N(v)\cup\{v\}$; in order to avoid ambiguity, we use also $N(v,G)$
instead of $N(v)$. For $A\subset V$, we denote $N(A)=\{v\in V-A:N(v)\cap
A\neq\varnothing\}$, and $N[A]=N(A)\cup A$. In particular, if $\left\vert
N(v)\right\vert =1$, then $v$ is a \textit{pendant vertex} of $G$, and
$pend(G)$ designates the set $\{v\in V(G):v$ \textit{is a pendant vertex in}
$G\}$. By \textit{tree} we mean a connected acyclic graph of order greater
than one, and a \textit{forest} is a disjoint union of trees and isolated
vertices. The bipartition $\{A,B\}$ of a tree $T$ is a partition of its set of
vertices into two stable sets $A$ and $B$. It is known that the bipartition
$\{A,B\}$ of a tree $T$ is unique up to isomorphism, and that $\alpha
(T)=\left\vert A\right\vert =\left\vert B\right\vert $ if and only if $T$ has
a perfect matching.

In this paper we show that any stable set $S$ of a tree $T$, of size
$\left\vert S\right\vert \geq\left\vert V(T)\right\vert /2$, contains at least
one pendant vertex of $T$. As a consequence, we infer that $S\cap
pend(T)\neq\varnothing$ is valid for any $S\in\Omega(T)$. Moreover, we prove
that in a tree $T$ with $\alpha(T)>\left\vert V(T)\right\vert /2$, there exist
at least two pendant vertices belonging to every maximum stable set of $T$,
such that the distance between them is even.

We give also a new proof for a result of Hopkins and Staton, \cite{HopStat},
stating that if $\{A,B\}$ is the standard bipartition of the vertex set of a
tree $T$, then $\Omega(T)=\{A\}$ or $\Omega(T)=\{B\}$ if and only if the
distance between any two pendant vertices of $T$ is even. Using this result we
deduce that if $\{A,B\}$ is the bipartition of a tree $T$ and $S$ is a stable
set such that $\left\vert S\right\vert >\min\{\left\vert A\right\vert
,\left\vert B\right\vert \}$, then $S\cap pend(T)\neq\varnothing$.

Our findings are also incorporated in the following contexts.

Firstly, the following theorem concerning maximum stable sets in general
graphs, due to Nemhauser and Trotter \cite{NemhTro}, shows that for a special
subgraph $H$ of a graph $G$, some maximum stable set of $H$ can be enlarged to
a maximum stable set of $G$. Namely, if $A\in\Omega(G[N[A]])$, then there is
$S\in\Omega(G)$, such that $A\subseteq S$. We show that, for trees, some kind
of an inverse theorem is also true. More precisely, we show that any maximum
stable set of a tree $T$ contains at least one of its pendant vertices, i.e.,
for any $S\in\Omega(T)$ there exists some $A\subseteq S$, such that $A$ is a
maximum stable set in the subgraph of $T$ induced by $N[A]$.

Secondly, in \cite{BorosGolLev} it was proved that if $G$ is connected, then
$\xi(G)\geq1+\alpha(G)-\mu(G)$, where $\mu(G)$ is the matching number of $G$.
This strengthened the following finding stated in \cite{LeviMan3}: if
$\alpha(G)>(\left\vert V(G)\right\vert +k-1)/2$, then $\xi(G)\geq k+1$;
moreover, $\xi(G)\geq k+2$ is valid, whenever $\left\vert V(G)\right\vert
+k-1$ is an even number. For $k=1$, the previous inequality provides us with a
generalization of a result of Hammer, Hansen and Simeone \cite{HamHanSim}
claiming that if a graph $G$ has $\alpha(G)>\left\vert V\left(  G\right)
\right\vert /2$, then $\xi(G)\geq1$. In \cite{LeviMan2} it was shown that if
$G$ is a connected bipartite graph with $\left\vert V(G)\right\vert \geq2 $,
then $\xi(G)\neq1$. Jamison \cite{Jamison}, Zito \cite{Zito}, and Gunther,
Hartnel and Rall \cite{GunHarRall} proved independently that $\xi(G)\neq1$ is
true for any tree $T$. Now, for a tree $T$ with $\alpha(T)>\left\vert
V(T)\right\vert /2$, we demonstrate that there exist at least two pendant
vertices even distance apart belonging to all maximum stable sets of $T$.

Thirdly, it is well-known that any tree $T$ has at least two pendant vertices
(see, for instance, Berge \cite{Berge}). Our results say that if
$\alpha(T)>\left\vert V(T)\right\vert /2$, then at least two pendant vertices
of $T$ belong to all maximum stable sets of $T$, and whenever $\alpha
(T)=\left\vert V(T)\right\vert /2$ then both parts of the bipartition of $T$
contain at least one pendant vertex.

\section{Pendant vertices and maximum stable sets}

\begin{lemma}
\label{lem2}Any stable set consisting of only pendant vertices of a graph $G$
is contained in a maximum stable set of $G$.
\end{lemma}

\begin{proof}
Let $A$ be a stable set of $G$ such that $A\subseteq pend(G)$, and $S\in
\Omega(G)$. If $u\in A-S$, then $u$ is adjacent to some $w\in S-A$, otherwise
$S\cup\left\{  u\right\}  $ is a stable set larger than $S$, which contradicts
the maximality of $S$. Hence, $S_{1}=S\cup\{u\}-\left\{  w\right\}  \in
\Omega(G)$, and $\left\vert A\cap S\right\vert <\left\vert A\cap
S_{1}\right\vert $. Therefore, using this exchange procedure, after a finite
number of steps, we have to obtain a maximum stable set including $A$.
\end{proof}

The converse of Lemma \ref{lem2} is not generally true. For instance, as it is
emphasized in Figure \ref{fig1}, the maximum stable set $S=\{a,b,c\}$ does not
contain any pendant vertex of the graph $G$.

\begin{figure}[h]
\setlength{\unitlength}{1.0cm} \begin{picture}(5,1.1)\thicklines
\put(5,0){\line(1,0){4}}
\multiput(5,0)(2,0){3}{\circle*{0.25}}
\multiput(6,0)(2,0){2}{\circle*{0.25}}
\put(6,1){\circle*{0.25}}
\put(5,0){\line(1,1){1}}
\put(6,0){\line(0,1){1}}
\put(7,1){\circle*{0.25}}
\put(7,0){\line(0,1){1}}
\put(6,1){\line(1,0){1}}
\put(6.2,0.3){\makebox(0,0){$a$}}
\put(7.3,1){\makebox(0,0){$b$}}
\put(8,0.3){\makebox(0,0){$c$}}
\put(4.5,0.5){\makebox(0,0){$G$}}
\end{picture}
\caption{$S=\{a,b,c\}\in\Omega(G)$ and $S\cap pend(G)=\varnothing$.}%
\label{fig1}%
\end{figure}

\begin{theorem}
\label{th4}If $S$ is a stable set of a tree $T=(V,E)$ and $\left\vert
S\right\vert \geq\left\vert V\right\vert /2$, then $S\cap pend(T)$ is not
empty. Moreover, if $S-pend(T)\neq\varnothing$, then there exist $v\in S\cap
pend(T)$ and $w\in S$, such that the distance between them equals two.
\end{theorem}

\begin{proof}
Suppose, on the contrary, that $S\cap pend(T)=\varnothing$. Hence, any $s\in
S$ has $\left\vert N(s)\right\vert \geq2$. Since $T$ is a tree and $S$ is a
stable set of size $\left\vert S\right\vert \geq\left\vert V\right\vert /2$,
it yields the following contradiction:
\[
\left\vert V\right\vert -1=\left\vert E\right\vert \geq\left\vert
(S,V-S)\right\vert \geq2\left\vert S\right\vert \geq\left\vert V\right\vert .
\]
Consequently, we infer that $S\cap pend(T)$ is not empty.

We can assert now that there exists some $k\geq1$, such that $S\cap
pend(T)=\{v_{i}:1\leq i\leq k\}$.

\textit{Case 1.} There exist two vertices from $S\cap pend(T)$ at distance two apart.

\textit{Case 2.} Any two vertices of $S\cap pend(T)$ are at distance at least three.

Let us denote $N(v_{i})=\{u_{i}\},1\leq i\leq k$. According to the hypothesis
of the case, all the vertices $u_{i}$ are different.

Assume, on the contrary, that for any $v_{i}\in S\cap pend(T)$ and any $w\in
S-pend(T)$, the distance between them is greater than two. Hence, any $w\in
S-pend(T)$ has $\left\vert N(w)\right\vert \geq2$.

Let $H_{1},H_{2},...,H_{p},1\leq p\leq k$, be the connected components of the
subgraph $T\left[  \bigcup\limits_{i=1}^{k}N\left[  v_{i}\right]  \right]  $,
of orders $n_{1},n_{2},...,n_{p}$, respectively. Since, by our assumption, no
$w\in S-pend(T)$ is connected to any $u_{i}$, we infer that for every $H_{j}$
there exists an edge joining this component to a vertex from $V-S$. Hence, it
yields the following contradiction:
\[
\left\vert V\right\vert -1=\left\vert E\right\vert \geq2\left(  \left\vert
S\right\vert -k\right)  +p+%
{\displaystyle\sum\limits_{i=1}^{p}}
\left(  n_{i}-1\right)  =2\left\vert S\right\vert \geq\left\vert V\right\vert
.
\]

Consequently, there must exist some $v\in S\cap pend(T)$ and $w\in S$ such
that the distance between them equals two.
\end{proof}

Let us notice that if the condition $S-pend(T)\neq\varnothing$ in Theorem
\ref{th4} is not satisfied, then all the distances between different vertices
of $S$ can be greater than two, e.g., see the stable set $\{v_{1},v_{2}%
,v_{3}\}$ of the tree $T_{1}$ in Figure \ref{fig323}.

On the other hand, if the condition $S-pend(T)\neq\varnothing$ in Theorem
\ref{th4} is satisfied, all the distances between the vertices of $S-pend(T)$
and vertices of $S$ can be different from $2$, for instance, see the stable
set $\{v_{1},v_{2},v_{3},v_{4},v_{5},v_{6},w_{1},w_{2}\}$ of the tree $T_{2}$
in Figure \ref{fig323}. 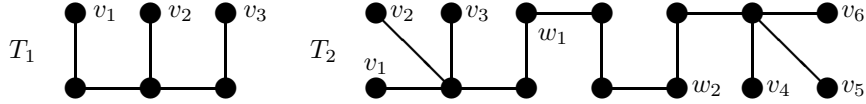
\begin{figure}[h]
\setlength{\unitlength}{1.0cm} \begin{picture}(5,1.2)\thicklines
\multiput(2,0)(1,0){3}{\circle*{0.29}}
\multiput(2,1)(1,0){3}{\circle*{0.29}}
\multiput(2,0)(1,0){3}{\line(0,1){1}}
\multiput(2,0)(1,0){2}{\line(1,0){1}}
\put(2.4,1){\makebox(0,0){$v_{1}$}}
\put(3.4,1){\makebox(0,0){$v_{2}$}}
\put(4.4,1){\makebox(0,0){$v_{3}$}}
\put(1.3,0.5){\makebox(0,0){$T_{1}$}}
\multiput(6,0)(1,0){7}{\circle*{0.29}}
\multiput(6,1)(1,0){7}{\circle*{0.29}}
\put(6,0){\line(1,0){2}}
\put(6,1){\line(1,-1){1}}
\put(7,0){\line(0,1){1}}
\put(8,0){\line(0,1){1}}
\put(8,1){\line(1,0){1}}
\put(9,0){\line(0,1){1}}
\put(9,0){\line(1,0){1}}
\put(10,0){\line(0,1){1}}
\put(10,1){\line(1,0){2}}
\put(11,0){\line(0,1){1}}
\put(11,1){\line(1,-1){1}}
\put(6,0.3){\makebox(0,0){$v_{1}$}}
\put(6.35,1){\makebox(0,0){$v_{2}$}}
\put(7.35,1){\makebox(0,0){$v_{3}$}}
\put(8.35,0.7){\makebox(0,0){$w_{1}$}}
\put(10.4,0){\makebox(0,0){$w_{2}$}}
\put(11.35,0){\makebox(0,0){$v_{4}$}}
\put(12.35,0){\makebox(0,0){$v_{5}$}}
\put(12.35,1){\makebox(0,0){$v_{6}$}}
\put(5.3,0.5){\makebox(0,0){$T_{2}$}}
\end{picture}
\caption{Distances between vertices belonging to maximal stable sets of
trees.}%
\label{fig323}%
\end{figure}

Now using the fact that $\alpha(G)\geq\left\vert V\left(  G\right)
\right\vert /2$ holds for any bipartite graph $G$, we obtain the following result.

\begin{corollary}
\label{cor8}If $T$ is a tree, then $S\cap pend(T)\neq\varnothing$ for any
$S\in\Omega(T)$.
\end{corollary}

\begin{figure}[h]
\setlength{\unitlength}{1.0cm} \begin{picture}(5,1.2)\thicklines
\multiput(4,0)(1,0){4}{\circle*{0.29}}
\multiput(4,1)(1,0){3}{\circle*{0.29}}
\multiput(4,0)(1,0){3}{\line(0,1){1}}
\multiput(4,0)(1,0){2}{\line(1,1){1}}
\multiput(6,0)(1,0){2}{\line(-1,1){1}}
\put(7,0){\line(-2,1){2}}
\put(3.7,0){\makebox(0,0){$a$}}
\put(4.7,0){\makebox(0,0){$b$}}
\put(6.3,0){\makebox(0,0){$c$}}
\put(7.3,0){\makebox(0,0){$d$}}
\put(3.3,0.5){\makebox(0,0){$G_{1}$}}
\multiput(9,0)(2,0){2}{\circle*{0.29}}
\put(10,0){\circle*{0.29}}
\put(9,0){\line(1,0){2}}
\multiput(9,1)(2,0){2}{\circle*{0.29}}
\put(10,1){\circle*{0.29}}
\put(10,1){\line(1,0){1}}
\multiput(9,0)(1,0){3}{\line(0,1){1}}
\put(8.7,0){\makebox(0,0){$u$}}
\put(9.7,1){\makebox(0,0){$v$}}
\put(11.3,0){\makebox(0,0){$w$}}
\put(8.3,0.5){\makebox(0,0){$G_{2}$}}
\end{picture}
\caption{Bipartite graphs with maximum stable sets containing no pendant
vertices.}%
\label{fig3}%
\end{figure}
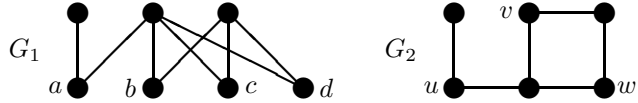

Corollary \ref{cor8} is not true for any connected graph $G$ with
$pend(G)\neq\varnothing$ (see, for instance, the graph in Figure \ref{fig1}).
Notice also that it cannot be generalized to a bipartite graph $G$ with
$pend(G)\neq\varnothing$, both for $\alpha(G)>\left\vert V(G)\right\vert /2 $
and $\alpha(G)=\left\vert V(G)\right\vert /2$ (see the graphs $G_{1},G_{2} $,
depicted in Figure \ref{fig3}, and the sets $\{a,b,c,d\}\in\Omega
(G_{1}),\{u,v,w\}\in\Omega(G_{2})$). 

\begin{corollary}
\label{cor4}If $T$ is a tree with $\alpha(T)=\left\vert V(T)\right\vert /2$,
then $T$ contains at least two pendant vertices at odd distance apart.
\end{corollary}

\begin{proof}
Let $\{A,B\}$ be a bipartition of $T$. Both $A$ and $B$ are maximum stable
sets, because $A,B$ are stable and $\alpha(T)\geq\max(\left\vert A\right\vert
,\left\vert B\right\vert )\geq\left\vert V\left(  T\right)  \right\vert /2$.
Hence, Corollary \ref{cor8} implies that both $A\cap pend(T)\neq\varnothing
\ $and$\ B\cap pend(T)\neq\varnothing$, that supports the conclusion.
\end{proof}

\begin{figure}[h]
\setlength{\unitlength}{1.0cm} \begin{picture}(5,1.2)\thicklines
\multiput(3,1)(1,0){5}{\circle*{0.25}}
\put(4,0){\circle*{0.25}}
\multiput(6,0)(1,0){3}{\circle*{0.25}}
\put(4,0){\line(-1,1){1}}
\put(4,0){\line(0,1){1}}
\put(4,0){\line(1,1){1}}
\put(4,0){\line(2,1){2}}
\put(6,0){\line(0,1){1}}
\put(6,0){\line(1,1){1}}
\put(7,0){\line(0,1){1}}
\put(8,0){\line(-1,1){1}}
\put(3.6,0){\makebox(0,0){$a$}}
\put(7.3,1){\makebox(0,0){$b$}}
\put(2.3,0.5){\makebox(0,0){$T_{1}$}}
\multiput(10,0)(1,0){2}{\circle*{0.25}}
\multiput(10,1)(1,0){3}{\circle*{0.25}}
\put(10,0){\line(0,1){1}}
\put(10,0){\line(1,1){1}}
\put(11,0){\line(0,1){1}}
\put(11,0){\line(1,1){1}}
\put(9.7,0){\makebox(0,0){$a$}}
\put(10.7,0){\makebox(0,0){$b$}}
\put(9.3,0.5){\makebox(0,0){$T_{2}$}}
\end{picture}
\caption{$S=\{a,b\}$ is a maximal stabe set containing no pendant vertices.}%
\label{fig3333}%
\end{figure}
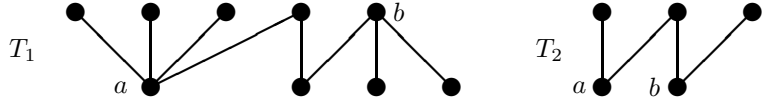

The tree $T_{1}$ in Figure \ref{fig3333} shows that the converse of Corollary
\ref{cor4} is not generally true.

Recall from \cite{HopStat} that $G$ is called a \textit{strong unique
independent graph} if $\left\vert \Omega(G)\right\vert =\left\vert
\{S\}\right\vert =1$ and $V(G)-S$ is also stable. For example, every cordless
path of odd order belongs to this class of graphs. Any strong unique
independent graph $G$ is necessarily bipartite, and its bipartition is
$\{S,V(G)-S\}$. Using Theorem \ref{th4}, we are giving now an alternative
proof of the following theorem characterizing strong unique independent trees,
which was first proved in \cite{HopStat}.

\begin{theorem}
\label{prop6}\cite{HopStat} If $\{A,B\}$ is the bipartition of the tree $T$,
then the following assertions are equivalent:

\emph{(i)} $T$ is a strong unique independent tree;

\emph{(ii)} $pend(T)\subseteq A$ or $pend(T)\subseteq B$;

\emph{(iii)} the distance between any two pendant vertices of $T$ is even.
\end{theorem}

\begin{proof}
\emph{(i)} $\Rightarrow$ \emph{(ii)} If $\Omega(T)=\{A\}$, then Lemma
\ref{lem2} implies $pend(T)\subseteq A$.

The equivalence \emph{(ii)} $\Leftrightarrow$ \emph{(iii)} is clear.

\emph{(ii)} $\Rightarrow$ \emph{(i)} Without loss of generality, we may
suppose that $\left\vert A\right\vert \geq\left\vert V(T)\right\vert /2$.
Since $A$ is also a stable set, Theorem \ref{th4} ensures that $A\cap
pend(T)\neq\varnothing$, and, consequently, $pend(T)\subseteq A$. Let
$S\in\Omega(T)$ and assume that $S-A\neq\varnothing$. Then
\[
\left\vert S\cap A\right\vert +\left\vert S\cap B\right\vert =\left\vert
S\right\vert \geq\left\vert A\right\vert =\left\vert S\cap A\right\vert
+\left\vert N(S\cap B)\right\vert
\]
and, hence, $\left\vert S\cap B\right\vert \geq\left\vert N(S\cap
B)\right\vert $. Since no vertex in $S\cap A$ is adjacent to any vertex in
$S\cap B$, it follows that at least one tree, say $T^{\prime}$, of the forest
$T[N[S\cap B]]$, has $\left\vert V(T^{\prime})\cap B\right\vert \geq\left\vert
V(T^{\prime})\cap A\right\vert $. Consequently, by Theorem \ref{th4}, it
proves that $T^{\prime}$ has at least one pendant vertex, say $v$, in
$V(T^{\prime})\cap B$. Since $N(v,T)=N(v,T^{\prime})$, we infer that some
pendant vertex of $T$ must be in $S\cap B$, in contradiction with
$pend(T)\subseteq A$.

Therefore, every maximum stable set of $T$ is a subset of $A$. Since $A$
is\emph{\ }stable, it follows that, in fact, $\Omega(T)=\{A\}$, i.e., $T$ is a
strong unique independent tree.
\end{proof}

\begin{theorem}
\label{th5} Let $\{A,B\}$ be the bipartition of the tree $T$. If $S$ is a
stable set such that $\left\vert S\right\vert >\min\{\left\vert A\right\vert
,\left\vert B\right\vert \}$, then $S\cap pend(T)\neq\varnothing$. Moreover,
there exist $v\in S\cap pend(T)$ and $w\in S$ such that the distance between
them is $2$.
\end{theorem}

\begin{proof}
We prove that if a stable set $S$ of $T$ satisfies $S\cap pend(T)=\varnothing
$, then $\min\{\left\vert A\right\vert ,\left\vert B\right\vert \}\geq
\left\vert S\right\vert $. If $S_{A}=S\cap A$ and $S_{B}=S\cap B$, then
$\{S_{A},B-S_{B}\}$ is the bipartition of the forest $F_{1}=T[S_{A}%
\cup(B-S_{B})]$, while $\{A-S_{A},S_{B}\}$ is the bipartition of the forest
$F_{2}=T[(A-S_{A})\cup S_{B}]$. Since $(S_{A},S_{B})=\varnothing$ and $S\cap
pend(T)=\varnothing$, it follows that $S_{A}\cap pend(F_{1})=\varnothing
=S_{B}\cap pend(F_{2})$. Consequently, by Theorem \ref{prop6}, every connected
component of $F_{1}$ or $F_{2}$, which is different from an isolated vertex,
is a strong unique independent tree. Moreover, every isolated vertex of
$F_{1}$ belongs to $B-S_{B}$, and every isolated vertex of $F_{2}$ belongs to
$A-S_{A}$. Therefore, both $\left\vert S_{B}\right\vert \leq\left\vert
A-S_{A}\right\vert $ and $\left\vert S_{A}\right\vert \leq\left\vert
B-S_{B}\right\vert $. Hence, we get that
\[
\left\vert S\right\vert =\left\vert S_{A}\right\vert +\left\vert
S_{B}\right\vert \leq\min\{\left\vert S_{A}\right\vert +\left\vert
A-S_{A}\right\vert ,\left\vert B-S_{B}\right\vert +\left\vert S_{B}\right\vert
\}=\min\{\left\vert A\right\vert ,\left\vert B\right\vert \},
\]
which completes the proof of the first assertion.

Now, let $S$ be a stable set such that
\[
\left\vert S\right\vert >\left\vert B\right\vert =\min\{\left\vert
A\right\vert ,\left\vert B\right\vert \},
\]
and let $S_{A}=S\cap A$, $S_{B}=S\cap B$. Any $v\in S_{A}$ has $\deg(v)\geq1$
and $N(v)\subseteq B-S_{B}$. Since
\[
\left\vert S_{A}\right\vert +\left\vert S_{B}\right\vert =\left\vert
S\right\vert >\left\vert B\right\vert =\left\vert S_{B}\right\vert +\left\vert
B-S_{B}\right\vert ,
\]
we see that $\left\vert S_{A}\right\vert >\left\vert B-S_{B}\right\vert $.
Therefore, it follows that some tree $H=(A_{1},B_{1},E_{1})$ of the forest
$T[S_{A}]$ must have $\left\vert A_{1}\right\vert >\left\vert B_{1}\right\vert
>1$, where $A_{1}\subseteq S_{A}$. According to Theorem \ref{th4}, it follows
that $A_{1}\cap pend(H)\neq\varnothing$, which implies that $A_{1}\cap
pend(T)\neq\varnothing$, because $A_{1}\cap pend(H)\subseteq A_{1}\cap
pend(T)$. In addition, if $A_{1}\subseteq pend(H)$, then there is
$\{v,w\}\subseteq A_{1}$ and the distance between $v,w$ equals $2$, while if
$A_{1}-pend(H)\neq\varnothing$, then, according to Theorem \ref{th4}, there
exist $v\in A_{1}\cap pend(H)$ and $w\in A_{1}$, such that the distance
between them is $2$. In both cases, we may conclude that there are $v\in S\cap
pend(T)$ and $w\in S$ such that the distance between them equals two.
\end{proof}

Theorem \ref{th5} shows that for a maximal stable set $S$ of a tree with the
bipartition $\{A,B\}$, it is enough to require that there are no pendant
vertices belonging to $S$ to ensure that $\min\{\left\vert A\right\vert
,\left\vert B\right\vert \}\geq\left\vert S\right\vert $. In\ Figure
\ref{fig3333} are depicted examples of trees with $\min\{\left\vert
A\right\vert ,\left\vert B\right\vert \}=\left\vert S\right\vert $ and
$\min\{\left\vert A\right\vert ,\left\vert B\right\vert \}>\left\vert
S\right\vert $. 

\begin{figure}[h]
\setlength{\unitlength}{1.0cm} \begin{picture}(5,1.2)\thicklines
\multiput(2,1)(1,0){5}{\circle*{0.25}}
\multiput(3,0)(2,0){2}{\circle*{0.25}}
\put(3,0){\line(-1,1){1}}
\put(3,0){\line(0,1){1}}
\put(3,0){\line(1,1){1}}
\put(5,0){\line(-1,1){1}}
\put(5,0){\line(0,1){1}}
\put(5,0){\line(1,1){1}}
\put(2.7,0){\makebox(0,0){$a$}}
\put(4.7,1){\makebox(0,0){$b$}}
\put(5.7,1){\makebox(0,0){$c$}}
\put(1.6,0.5){\makebox(0,0){$T_{1}$}}
\multiput(7.5,1)(1,0){5}{\circle*{0.25}}
\multiput(8.5,0)(1,0){5}{\circle*{0.25}}
\put(8.5,0){\line(-1,1){1}}
\put(8.5,0){\line(0,1){1}}
\put(8.5,0){\line(1,1){1}}
\put(9.5,0){\line(0,1){1}}
\put(9.5,0){\line(1,1){1}}
\put(10.5,0){\line(0,1){1}}
\put(11.5,0){\line(-1,1){1}}
\put(11.5,0){\line(0,1){1}}
\put(11.5,1){\line(1,-1){1}}
\put(8.8,0){\makebox(0,0){$a$}}
\put(10.8,1){\makebox(0,0){$b$}}
\put(12.8,0){\makebox(0,0){$c$}}
\put(7.1,0.5){\makebox(0,0){$T_{2}$}}
\end{picture}
\caption{$S=\{a,b,c\}$ is a maximal stabe set in $T_{1},T_{2}$, containing
pendant vertices.}%
\label{fig3434}%
\end{figure}
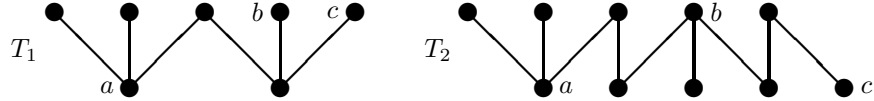

If $\left\vert A\right\vert \neq\left\vert B\right\vert $ then the claim of
Theorem \ref{th5} is stronger than the corresponding direct consequence from
Theorem \ref{th4}, because there is a tree $T$ containing a maximal stable set
$S$, such that $\min\{\left\vert A\right\vert ,\left\vert B\right\vert
\}<\left\vert S\right\vert <\left\vert V(T)\right\vert /2$ and $S\cap
pend(T)\neq\varnothing$ (for an example, see the tree $T_{1}$ in Figure
\ref{fig3434}). 

The tree $T_{2}$ in Figure \ref{fig3434} shows that the converse of Theorem
\ref{th5} is not true.

\section{Pendant vertices and intersection of all maximum stable sets}

Recall the following result, which we shall use in the sequel.

\begin{proposition}
\cite{LeviMan3}\label{prop4} For a connected bipartite graph $G$ of order at
least two, the following assertions are true:

\emph{(i)} $\alpha(G)>\left\vert V\left(  G\right)  \right\vert /2$ if and
only if $\xi(G)\geq2$;

\emph{(ii)} $\alpha(G)=\left\vert V\left(  G\right)  \right\vert /2$ if and
only if $\xi(G)=0$.
\end{proposition}

Let $G_{i}=(V_{i},E_{i}),i=1,2$, be two graphs with $V_{1}\cap V_{2}%
=\varnothing$, and $Q_{1},Q_{2}$ be cliques of the same size in $G_{1},G_{2}$,
respectively. The \textit{clique bonding} of the graphs $G_{1},G_{2}$ is the
graph $G=G_{1}\ast Q\ast G_{2}$ obtained by identifying $Q_{1}$ and $Q_{2}$
into a single clique $Q$, \cite{Berge}. If $V(Q_{1})=\{v_{1}\}$,
$V(Q_{2})=\{v_{2}\}$, we shall denote the clique bonding of $G_{1}$ and
$G_{2}$ by $G_{1}\ast v\ast G_{2}$. In other words, $V(G_{1}\ast v\ast
G_{2})=V_{1}\cup V_{2}\cup\left\{  v\right\}  -\left\{  v_{1},v_{2}\right\}  $
and
\[
E(G_{1}\ast v\ast G_{2})=E\left(  G_{1}\left[  V_{1}-v_{1}\right]  \right)
\cup E\left(  G_{2}\left[  V_{2}-v_{2}\right]  \right)  \cup\{vx:v_{1}x\in
E_{1}\}\cup\{vy:v_{2}y\in E_{2}).
\]

\begin{lemma}
\label{lem3}Let $T_{1},T_{2}$ be trees and $T=T_{1}*v*T_{2}$.

\emph{(i)} if $v\in core(T)$, then $\alpha(T)=\alpha(T_{1})+\alpha(T_{2})-1$;

\emph{(ii)} $v\in core(T)$ if and only if $v\in core(T_{i}),$ $i=1,2$;

\emph{(iii)} if $v\in core(T)$, then $core(T)=core(T_{1})\cup core(T_{2})$.
\end{lemma}

\begin{proof}
\emph{(i) }Let $S\in\Omega(T)$. Then $S\cap V(T_{i})$ is stable in $T_{i}$,
and, therefore, it follows that $\left\vert S\cap V(T_{i})\right\vert
\leq\alpha(T_{i})$, for each $i=1,2$. Hence, we get that
\begin{align*}
\alpha(T)  & =\left\vert S\right\vert =\left\vert S\cap V(T_{1})\right\vert
+\left\vert S\cap V(T_{2})-\{v\}\right\vert =\\
& =\left\vert S\cap V(T_{1})\right\vert +\left\vert S\cap V(T_{2})\right\vert
-1\leq\alpha(T_{1})+\alpha(T_{2})-1.
\end{align*}

\textit{Case 1}. There are $S_{i}\in\Omega(T_{i}),i=1,2$, such that $v\in
S_{1}\cap S_{2}$. Then $S_{1}\cup S_{2}$ is stable in $T$ and
\[
\left\vert S_{1}\cup S_{2}\right\vert =\alpha(T_{1})+\alpha(T_{2})-1\geq
\alpha(T),
\]
and this implies that $\alpha(T)=\alpha(T_{1})+\alpha(T_{2})-1$.

\textit{Case 2}. There are $S_{i}\in\Omega(T_{i})$, such that $v\notin
S_{1}\cap S_{2}$, i.e., $v\notin S_{1}$, (or $v\notin S_{2}$). Hence,
$S_{3}=S_{1}\cup S_{2}-\{v\}$ is stable in $T$, and
\[
\left\vert S_{3}\right\vert \geq\alpha(T_{1})+\alpha(T_{2})-1\geq\alpha(T),
\]
and this leads to the following contradiction with the hypothesis $v\in
core(T)$:\emph{\ }$S_{3}\in\Omega(T)$ and $v\notin S_{3}$.

Thus, we may conclude that $\alpha(T)=\alpha(T_{1})+\alpha(T_{2})-1$.

\emph{(ii) }If $v\in core(T)$, then \textit{Case 2} of part \emph{(i)}
explicitly means that $v\in core(T_{i}),$ $i=1,2$.

Conversely, let $v\in core(T_{i}),$ $S_{i}\in\Omega(T_{i}),i=1,2$, and assume
that there is $S\in\Omega(T)$, such that $v\notin S$. Then, the set $S_{1}\cup
S_{2}$ is stable in $T$ and $\left\vert S_{1}\cup S_{2}\right\vert
=\alpha(T_{1})+\alpha(T_{2})-1$. Clearly, $S\cap V(T_{i})$ is stable in
$T_{i},i=1,2$, and because $v\in core(T_{i}),i=1,2$, we have that $\left\vert
S\cap V(T_{i})\right\vert \leq\alpha(T_{i})-1,i=1,2$. Hence,
\[
\left\vert S\right\vert =\left\vert S\cap V(T_{1})\right\vert +\left\vert
S\cap V(T_{2})\right\vert \leq\alpha(T_{1})+\alpha(T_{2})-2<\left\vert
S_{1}\cup S_{2}\right\vert ,
\]
and this contradicts the choice $S\in\Omega(T)$.

\emph{(iii) }According to part \emph{(i)}, we have $\alpha(T)=\alpha
(T_{1})+\alpha(T_{2})-1$, and part \emph{(ii)} ensures that $v\in
core(T_{i}),i=1,2$.

Let $w\in(core(T)-\{v\})\cap V(T_{1})$ and $S_{i}\in\Omega(T_{i}),i=1,2$. Then
$S_{1}\cup S_{2}\in\Omega(T)$, and, therefore, $w\in S_{1}$. Since $S_{1}$ is
an arbitrary set from $\Omega(T_{1})$, we get that $w\in core(T_{1}).$

Similarly, one can show that if $w\in(core(T)-\{v\})\cap V(T_{2})$, then $w\in
core(T_{2})$. Therefore, we may conclude that $core(T)\subseteq core(T_{1}%
)\cup core(T_{2}).$

Conversely, let $w\in core(T_{1})-\{v\}$, and suppose there is $S\in\Omega
(T)$, such that $w\notin S$. Let us denote $S_{i}=S\cap V(T_{i})$, for $i=1,2
$. Since $w\notin S_{1}$, it follows that $\left\vert S_{1}-\{v\}\right\vert
\leq\alpha(T_{1})-2$. Hence, we get a contradiction:
\[
\left\vert S\right\vert =\left\vert S_{1}-\{v\}\right\vert +\left\vert
S_{2}\right\vert \leq\alpha(T_{1})-2+\alpha(T_{2})<\alpha(T_{1})+\alpha
(T_{2})-1=\alpha(T)=\left\vert S\right\vert .
\]
Consequently, $core(T_{1})\cup core(T_{2})\subseteq core(T)$ is also valid,
and this completes the proof.
\end{proof}

In the following statement we are strengthening Proposition \ref{prop4} for
the case of trees.

\begin{theorem}
\label{th6}If $T$ is a tree with $\alpha(T)>\left\vert V(T)\right\vert /2$,
then $\left\vert core(T)\cap pend(T)\right\vert \geq2$.
\end{theorem}

\begin{proof}
According to Proposition \ref{prop4}\emph{(i)}, we infer that $\xi(T)\geq2$.
Since $T$ is a tree, it follows that $\left\vert pend(T)\right\vert \geq2$.

To prove the theorem we use induction on $n=\left\vert V(T)\right\vert $. The
result is clearly true for $n=3$.

Let $T=(V,E)$ be a tree with $n=\left\vert V\right\vert >3$, and suppose that
the assertion is valid for any tree with fewer vertices.\ If $core(T)=pend(T)$%
, the result is clear. If $core(T)\neq pend(T)$, let $v\in core(T)-pend(T)$
and $T_{1},T_{2}$ be two trees such that $T=T_{1}\ast v\ast T_{2}$. A
partition of $N(v)$ in two non-empty sets gives rise to a corresponding
division of $T$ into $T_{1}$ and $T_{2}$.\emph{\ }According to Lemma
\ref{lem3}\emph{(ii)}, $v\in core(T_{i}),i=1,2$. Hence, Proposition
\ref{prop4}\emph{(i)} implies that $\alpha(T_{i})>\left\vert V(T_{i}%
)\right\vert /2,i=1,2$. By the induction hypothesis, each $T_{i}$ has at least
two pendant vertices belonging to $core(T_{i})$. Lemma \ref{lem3}\emph{(iii)}
ensures that $core(T)=core(T_{1})\cup core(T_{2})$, and, therefore, $T$ itself
has at least two pendant vertices in $core(T)$.
\end{proof}

\begin{corollary}
Let $T$ be a tree with $\alpha(T)>\left\vert V\left(  T\right)  \right\vert /2
$, and $k\geq2$. If there is a vertex $v\in core(T)$ of $deg(v)$ $\geq2k $,
then $\left\vert core(T)\cap pend(T)\right\vert \geq2k$.
\end{corollary}

\begin{proof}
Let us partition $N(v)$ into $k$ subsets $N_{i}(v),1\leq i\leq k$, each one
having at least two vertices. Then we can write
\[
T=(...((T_{1}\ast v\ast T_{2})\ast v\ast T_{3})...)\ast v\ast T_{k},
\]
where $T_{i}$ is the subtree of $T$ containing $N_{i}(v)$ as the neighborhood
of $v$.

Since $v$ is pendant in no $T_{i}$, we get $pend(T)=\cup\{pend(T_{i}):1\leq
i\leq k\}$.

By Lemma \ref{lem3}\emph{(iii)}, it follows that $core(T)=\cup\{core(T_{i}%
):1\leq i\leq k\}$.

According to Lemma \ref{lem3}\emph{(ii)}, the vertex $v\in core(T_{i}),1\leq
i\leq k$, and, consequently, Theorem \ref{th6} implies that:
\[
\left\vert core(T)\cap pend(T)\right\vert =\left\vert core(T_{1})\cap
pend(T_{1})\right\vert +...+\left\vert core(T_{k})\cap pend(T_{k})\right\vert
\geq2k,
\]
and this completes the proof.
\end{proof}

Let us remark that for every natural number $k$ there exists a tree $T$ with a
vertex $v$ of degree $k$ such that $v\in core(T)$. For instance, such a tree
$T=(V,E)$ can be defined as follows: $V=\{v\}\cup\{x_{i}:1\leq i\leq2k\}$ and
$E=\{vx_{i}:1\leq i\leq k\}\cup\{x_{i}x_{i+k}:1\leq i\leq k\}$.

\begin{theorem}
If $T$ is a tree with $\alpha(T)>\left\vert V(T)\right\vert /2$, then for at
least two distinct vertices from $core(T)\cap pend(T)$ the distance between
them is even. Moreover, if the set $core(T)\cap pend(T)$ contains exactly two
vertices, then the distance between them never equals four.
\end{theorem}

\begin{proof}
Let $\{A,B\}$ be the bipartition of $T$. Notice that the distance between two
vertices is even if and only if both of them belong to one set of the bipartition.

To prove the theorem we use induction on $n=\left\vert V(T)\right\vert $.

If $n=3$, then $T=P_{3}$ and the assertion is true.

Let now $T$ be a tree with $n\geq4$ vertices. According to Theorem \ref{th6},
$\alpha(T)>n/2$ yields $\left\vert core(T)\cap pend(T)\right\vert \geq2$.

\textit{Case 1.} $\left\vert core(T)\cap pend(T)\right\vert \geq3$.

Then we get%
\[
\min(\left\vert A\cap core(T)\cap pend(T)\right\vert ,\left\vert B\cap
core(T)\cap pend(T)\right\vert )>1.
\]
Hence, at least two vertices of $core(T)\cap pend(T)$ belong to one set of the
bipartition, i.e., the distance between them is even.

\textit{Case 2.} $\left\vert core(T)\cap pend(T)\right\vert =\left\vert
\{u,v\}\right\vert =2$. Figure \ref{fig4444} shows that such trees exist.
\begin{figure}[h]
\setlength{\unitlength}{1.0cm} \begin{picture}(5,1.4)\thicklines
\multiput(4,1)(1,0){3}{\circle*{0.25}}
\multiput(4,0)(1,0){6}{\circle*{0.25}}
\put(4,1){\line(1,0){2}}
\put(4,0){\line(1,0){5}}
\put(6,1){\line(1,-1){1}}
\put(3.7,1){\makebox(0,0){$u$}}
\put(3.7,0){\makebox(0,0){$v$}}
\put(5,1.3){\makebox(0,0){$x$}}
\put(5,0.35){\makebox(0,0){$y$}}
\put(6,1.3){\makebox(0,0){$c$}}
\put(6,0.4){\makebox(0,0){$d$}}
\end{picture}
\caption{A tree $T$ with $\left\vert core(T)\cap pend(T)\right\vert =2$.}%
\label{fig4444}%
\end{figure}

If $N(u)=N(v)$, then the distance between them is two, which is both even and
different from four.

Suppose now that $N(u)=\{x\}\neq\{y\}=N(v)$, and let $F=T[A\cup
B-\{u,v,x,y\}]$. Since $u$ and $v$ belong to all maximum stable sets of $T$,
we conclude that neither $x$ nor $y$ are contained in any maximum stable set
of $T$. Hence $\Omega(T)=\{S\cup\{u,v\}:S\in\Omega(F)\}$. Consequently,
$core(T)=core(F)\cup\{u,v\}$ and $\alpha(F)=\alpha
(T)-2>n/2-2=(n-4)/2=\left\vert V(F)\right\vert /2$. Suppose that $F$ consists
of $k\geq1$ trees $\{T_{i}:1\leq i\leq k\}$. Since $\alpha(F)=\alpha
(T_{1})+...+\alpha(T_{k})>\left\vert V(F)\right\vert /2$, at least one tree,
say $T_{j}$, has $\alpha(T_{j})>\left\vert V(T_{j})\right\vert /2$. By the
induction hypothesis, there exist two distinct vertices $c,d\in core(T_{j}%
)\cap pend(T_{j})$ such that the distance between them in $T_{j}$ is even.

We claim that $\{c,d\}\subset N(x)\cup N(y)$. Otherwise, if, for instance,
$c\notin N(x)\cup N(y)$, then $c\in core(T)\cap pend(T)$ and this contradicts
the fact that $core(T)\cap pend(T)=\{u,v\}$. Further, if $\{c,d\}\subset N(x)
$ or $\{c,d\}\subset N(y)$, then $T_{j}$ is not a tree, since $\{cx,xd\} $ or
$\{cy,yd\}$, respectively, builds a new path connecting $c$ and $d$ in
addition to the unique path between $c$ and $d$ in $T_{j}$ (together the two
paths create a cycle, which is forbidden in trees). Suppose that $c\in N(x)$
and $d\in N(y)$. Then $xy\notin$ $E(T)$, because, otherwise, $T_{j}$ can not
be a tree.

No edge from the set $\{uv,uy,ud,uc,vx,vc,vd\}$ exists in $T$, since the
vertices $u$ and $v$ are pendant in $T$. One can find an example of such a
situation in Figure \ref{fig4444}. The vertices $c$ and $d$ are not adjacent
in $T$, because they are pendant in $T_{j}$. Therefore, $\left\vert N(x)\cup
N(y)\right\vert \geq4$, and the shortest path between $u$ and $v$ goes through
the vertices $x,c,y,d$, at least. Thus, the distance between $u,v$ in $T$ is
greater than the distance between $c$ and $d$ in $T_{j}$ by four, and
consequently, it is even and, moreover, different from four, because the
vertices $c$ and $d$ are distinct.
\end{proof}

\section{Conclusions}

In this paper we have studied relationships between pendant vertices and
maximum stable sets of a tree. We have obtained a more precise version of the
well-known result of Berge, \cite{Berge}, stating that $\left\vert
pend(T)\right\vert \geq2$ holds for any tree $T$ having at least two vertices.
Namely, we have proved that for such a tree $T$ either it has a perfect
matching, and then both $A\cap pend(T)\neq\varnothing$ and $B\cap
pend(T)\neq\varnothing$, where $\{A,B\}$ is its bipartition, or it has no
perfect matching, and then at least two of its pendant vertices an even
distance apart belong to all maximum stable sets.

As open problems, we suggest the following.

Suppose that the tree $T$ has no perfect matching. Are there at least two
pendant vertices of $T$ belonging to $\cap\{S:S$ \textit{is a maximal stable
set in }$T$\textit{\ of size} $k\}$, for either $k=\left\vert V(T)\right\vert
/2$, or $k=\min\{\left\vert A\right\vert ,\left\vert B\right\vert \}$?

\end{document}